\documentclass[12pt]{article}
\usepackage{amsmath}
\usepackage{amssymb}
\usepackage[cp1251]{inputenc}
\usepackage[russian]{babel}
\newcommand{\il}[2]{\int\limits_{#1}^{#2}}
\newcommand{\ilp}[1]{\int\limits_{#1}^{+\infty}}

\newcommand{\ph}{\phantom{a}}
\newcommand{\phh}{\phantom{aaa}}

\newcommand{\sist}[2]{\left\{
\begin{array}{l}
{#1}\\
\ph\\
{#2}
\end{array}
\right.}

\begin{document}

MSC 34C10

\vskip 20pt

\centerline{\bf Kamenev type conditions for oscillation of third order }
\centerline{\bf linear ordinary differential equations}

\vskip 10 pt

\centerline{\bf G. A. Grigorian}

\vskip 10 pt

\centerline{0019 Armenia c. Yerevan, str. M. Bagramian 24/5}
\centerline{Institute of Mathematics of NAS of Armenia}
\centerline{E - mail: mathphys2@instmath.sci.am, \ph phone: 098 62 03 05, \ph 010 35 48 61}

\vskip 20 pt

\noindent
Abstract. The  Riccati equation method  is used to establish Kamenev-type conditions for the existence of oscillatory solutions to third order linear ordinary differential equations. Three  oscillatory theorems are proved, which generalize the Lazer's oscillation criterion.

\vskip 20 pt


\vskip 20 pt

{\bf 1. Introduction}. Let $p(t),\ph q(t)$ and $r(t)$ be real-valued continuous functions on $[t_0,+\infty)$. Consider the third order linear ordinary differential equation
$$
\phi''' + p(t)\phi'' + q(t)\phi' + r(t)\phi = 0, \phh t\ge t_0.
$$
Throughout this paper we will assume that $p(t)$ is  continuously differentiable on $[t_0,+\infty)$.

\vskip10pt

{\bf Definition 1.1.} {\it A solution of Eq. (1.1) is called oscillatory if it has arbitrarily large zeroes on $[t_0,+\infty)$.
}
\vskip10pt

The study of the oscillatory behavior of solutions to linear third-order ordinary \linebreak differential equations is an important problem in the qualitative theory of differential equations, and many works are devoted to it (see [1,2,4] and the works cited therein).  Among them notice the following result due to Lazer.

{\bf Theorem 1.1([4, Theorem 1.3]).} {\it If $p(t)\equiv 0, \ph q(t)\le 0, \ph r(t) > 0, \ph t \ge t_0$ and
$$
\ilp{t_0}\Bigl[r(t) - \frac{2}{3\sqrt{3}}(-q(t))^{3/2}\Bigr]d t = +\infty,
$$
then Eq. (1.1) has oscillatory solutions.
}
\phantom{aaaaaaaaaaaaaaaaaaaaaaaaaaaaaaaaaaaaaaa}$\blacksquare$

In this paper, we use the Riccati equation method to establish Kamenev-type conditions for the existence of oscillatory solutions to the equation. (1.1). Three oscillatory theorems generalizing Theorem 1.1 are proved.

{\bf Auxiliary propositions}.  The following two lemmas are of fundamental importance for the proof of the main results.

{\bf Lemma 2.1.([1, Lemma 2.2])} {\it If $q(t) \le 0, \ph r(t) > 0, \ph t \ge t_0$ and $\phi(t)\not\equiv 0$ is a solution of Eq. (1.1) with $\phi(t)\ge 0$ or $\phi(t)\le 0$ eventually, then there exists $T\in [t_0,+\infty)$ such that either
$$
\phi(t)\phi'(t)\le 0, \phh t \ge T \eqno (2.1)
$$
or
$$
\phi(t)\phi'(t) \ge 0, \ph t \ge T \ph \mbox{and} \ph \phi(t) > 0 \ph  \mbox{for} \ph t \ge T.
$$
Furthermore if (2.1) holds, then
$$
\phi(t)\phi'(t)\phi''(t) \ne 0, \ph sgn \ph \phi(t) = sgn\ph  \phi''(t) \ne sgn \ph \phi'(t), \phh t \ge t_0 \eqno (2.2)
$$
and
$$
\lim\limits_{t \to +\infty}\phi'(t) = \lim\limits_{t \to +\infty}\phi''(t) = 0, \ph \lim\limits_{t \to +\infty}\phi(t) = k \ne\pm\infty. \eqno (2.3)
$$
}

\phantom{aaaaaaaaaaaaaaaaaaaaaaaaaaaaaaaaaaaaaaaaaaaaaaaaaaaaaaaaaaaaaaaaaaaaaaaa}$\blacksquare$

{\bf Lemma 2.2([1. Lemma 2.3]).} {\it Let $q(t) \le 0, \ph r(t) > 0, \ph t \ge t_0.$ A necessary and sufficient condition for Eq. (1.1) to have oscillatory solutions is that for any nontrivial solution (2.2) and (2.3) hold.
}

\phantom{aaaaaaaaaaaaaaaaaaaaaaaaaaaaaaaaaaaaaaaaaaaaaaaaaaaaaaaaaaaaaaaaaaaaaaaa}$\blacksquare$

Let $f_k(t), \ph g_k(t), \ph h_k(t), \ph k=1,2,$ be real-valued continuous functions on $[t_0,+\infty)$. Consider the Riccati equations
$$
y' + f_k(t) y^2 + g_k(t) y + h_k(t) = 0, \phh t\ge t_0, \eqno (2.4_k)
$$
$k=1,2$ and the differential inequalities
$$
\eta' + f_k(t) \eta^2 + g_k(t) \eta + h_k(t) \ge 0, \phh t\ge t_0, \eqno (2.5_k)
$$
$k=1,2$.

{\bf Remark 2.1.} {\it Every solution  of Eq. $(2.4_2)$ on $[t_0,t_1)$ is also a solution of the inequality $(2.5_2)$ on $[t_0,t_1)$.}

{\bf Remark 2.2.} {\it If $f_1(t) \ge 0, \ph t\in [t_0,t_1)$, then every solution of the linear equation
$$
\zeta' + g_1(t)\zeta + h_1(t) = 0, \phh t\in [t_0,t_1)
$$
is also a solution of the inequality $(2.5_1)$ on $[t_0,t_1)$.}

{\bf Theorem 2.1 [3, Theorem 3.1].} {\it  Let $y_2(t)$ be a solution of Eq. $(2.4_2)$ on $[t_0,\tau_0) \linebreak (t_0 < \tau_0 \le +\infty)$ and let $\eta_1(t)$ and $\eta_2(t)$ be solutions of the inequalities $(2.5_1)$ and $(2.5_2)$  respectively on $[t_0,\tau_0)$ such that $y_2(t_0) \le \eta_k(t_0) \ph k=1,2.$ In addition let the following conditions be satisfied: $f_1(t) \ge 0, \ph \gamma - y_2(t_0) + \il{t_0}{t}\exp\biggl\{\il{t_0}{\tau}[f_1(s)(\eta_1(s) + \eta_2(s)) + g_1(s)]ds\biggr\}\biggl[(f_2(\tau) - f_1(\tau))^2 y_2^2(\tau) + (g_2(\tau) - g_1(\tau)) y_2(\tau) + h_2(\tau) - h_1(\tau)\biggr] d \tau \ge 0, \ph t\in [t_0,\tau_0)$ for some $\gamma \in [y_2(t_0), \eta_1(t_0)]$. Then Eq. $(2.4_1)$ has a solution $y_1(t)$ on $[t_0,\tau_0)$ with $y_1(t_0) \ge \gamma$ and $y_1(t) \ge y_2(t), \ph t\in [t_0,\tau_0)$.}

\phantom{aaaaaaaaaaaaaaaaaaaaaaaaaaaaaaaaaaaaaaaaaaaaaaaaaaaaaaaaaa}$\blacksquare$

{\bf 3. Oscillation criteria}. We set:

$$p_{-}(t)\equiv \min\{p(t),0\}, \ph t \ge t_0,\ph F(t,u) \equiv p(t) u^2 + (q(t) - p'(t)) u + r(t), \ph t \ge t_0, \ph u\ge 0$$
One can easily verify that for every fixed $t \ge t_0$ the minimum of the function $F$ is
$$
D(t)\equiv \sist{r(t), \ph \mbox{if} \ph p^2(t) + 3(p'(t) - q(t)) < 0,}{F_{min}(t), \ph \mbox{if}\ph p^2(t) + 3(p'(t) - q(t)) \ge 0, \ph t \ge t_0,}
$$
where $F_{min}(t)\equiv \min\Bigl\{r(t), \Bigl[\frac{\sqrt{p^2(t) + 3(p'(t) - q(t))} - p(t)}{3}\Bigr]^3 + p(t)\Bigl[\frac{\sqrt{p^2(t) + 3(p'(t) - q(t))} - p(t)}{3}\Bigr]^2 +\linebreak + (q(t) - p'(t)) \frac{\sqrt{p^2(t) +  3(p'(t) - q(t))} - p(t)}{3} + r(t)\bigr\},\ph t \ge t_0.$

In the case $p(t)\equiv 0, \ph q(t)\le 0, \ph r(t) > 0$ we have
$$
D(t) = r(t) - \frac{2}{3\sqrt{3}}(-q(t))^{3/2}, \phh t \ge t_0. \eqno (3.1)
$$

{\bf Theorem 3.1.} {\it Let the following conditions be satisfied.

\noindent
$(a) \ph r(t) > 0, \ph q(t) \le 0, \ph t \ge t_0$.

\noindent
$$
(b)\phantom{aaa} \limsup\limits_{t \to +\infty}\frac{1}{t^{\alpha + 1}}\il{t_0}{t}(t -\tau)^\alpha\biggl[(t - \tau) D(\tau) - \frac{\alpha+1}{9}p_{-}^2(\tau)\biggr]d \tau = +\infty, \ph \mbox{for some} \ph  \alpha > 0. \phantom{aaaaaaaaaaaaaaaaaaaaa}
$$
Then Eq. (1.1) has oscillatory solutions.
}

Proof. By Lemma 2.1 it follows from the conditions $(a)$ that for every nontrivial solution $\phi(t)$ of Eq. (1.1) with $\phi(t) \ge 0$ or $\phi(t) \le 0$ eventually there exists $T\ge t_0$ such that  either
$$
\phi'(t)\phi(t) \le 0, \phh t \ge T
$$
or
$$
\phi'(t)\phi(t) \ge 0, \ph t \ge T \ph \mbox{and}\ph \phi(t) > 0 \ph \mbox{for}\ph t\ge T. \eqno (3.2)
$$
By virtue of Lemma 2.2 to prove the theorem it is enough to show that the relations (3.2) are impossible. Suppose, at contrary, that (3.2) is valid. Then $y(t)\equiv \frac{\phi'(t)}{\phi(t)}, \ph t \ge T$ is a solution to the second order Riccati equation
$$
y'' + 3y[p(t) + y] y' + y^3 + p(t) y^2 + q(t) y + r(t) = 0, \phh t \ge T.
$$
Therefore, for $y(t)$ the following equality holds.
$$
y'(t) + \frac{3}{2}y^2(t) + p(t) y(t) + \il{T}{t}[y^3(\tau) + p(\tau) y^2(\tau) + (q(\tau) - p'(\tau)) y(\tau) + r(\tau)] d \tau = c_T, \ph t \ge T,
$$
where $c_T\equiv y'(T) +\frac{3}{2}y^2(T) + p(t) y(T)$.  Since $y(t) \ge 0,\ph t \ge T$ from the last equality we obtain $$
y'(t) + \frac{3}{2}y^2(t) + p_{-}(t) y(t) \le - \il{T}{t}D(\tau) d\tau + c_{T}, \phh t \ge T. \eqno (3.3).
$$
This inequality we can rewrite in the form
$$
y'(t) + \frac{3}{2}\Bigl(y(t) + \frac{p_{-}(t)}{3}\Bigr)^2 \le - \il{T}{t}D(\tau) d\tau + \frac{p^2_{-}(t)}{9}   + c_{T}, \phh t \ge T.
$$
Therefore,
$$
y(t) \le - \il{T}{t}d\tau\il{T}{\tau} D(\xi) d \xi + \frac{1}{9}\il{T}{t} p_{-}^2(\tau) d \tau + c_T(t - T), \ph t \ge T. \eqno (3.4)
$$
Without loss of generality we can take that $T > 0$. For any $\alpha > 0$ and $M > 0$ consider the integral operator
$$
K_{M,\alpha} \phi(t) \equiv \frac{\alpha(\alpha + 1)}{t^{\alpha+1}}\il{T}{t}(t -\tau)^{\alpha -1}\phi(\tau) d\tau, \ph t \ge M \ph \phi \in \mathbb{C}([M,+\infty)).
$$
Obviously this operator is monotone in the sense that if $\phi_j(t) \in \mathbb{C}([M,+\infty)), \ph j=~1,2, \linebreak \phi_1(t) \ge \phi_2(t), \ph t \ge M$, then $(K_{M,\alpha}\phi_1)(t) \ge (K_{M,\alpha}\phi_2)(t), \ph  t \ge M.$ Due to this, acting on both sides of (3.4) by operator $K_{T,\alpha}$ and making some simplifications we obtain
$$
\frac{\alpha(\alpha +1)}{t^{\alpha +1}}\il{T}{t}(t - \tau)^{\alpha -1} y(\tau) d \tau \le - \frac{1}{t^{\alpha+1}}\il{T}{t}(t - \tau)^{\alpha +1} D(\tau) d\tau + \frac{(\alpha +1)}{9 t^{\alpha +1}}\il{T}{t}(t - \tau)^{\alpha} p_{-}^2(\tau) d \tau +
$$
$$
+ \frac{\alpha(\alpha +1)c_T}{t^{\alpha +1}}\il{T}{t}(t - \tau)^{\alpha -1}(\tau - T) d \tau  + \frac{\alpha(\alpha +1)}{t^{\alpha +1}}y(T)\il{T}{t}(t - \tau)^{\alpha -1}  d \tau, \phh t \ge T. \eqno (3.5)
$$
Consider the function
$$
\Delta(t) \equiv \frac{1}{t^{\alpha+1}}\il{t_0}{T}(t- \tau)^{\alpha +1} D(\tau) d \tau - \frac{\alpha +1}{9 t^{\alpha +1}} \il{t_0}{T} (t - \tau)^\alpha p_{-}^2(\tau) d \tau +
$$
$$
+ \frac{\alpha(\alpha +1)c_T}{t^{\alpha +1}}\il{T}{t}(t -\tau)^{\alpha -1}(\tau - T) d\tau  + \frac{\alpha(\alpha +1)}{t^{\alpha +1}}y(T)\il{T}{t}(t - \tau)^{\alpha}  d \tau \ph t \ge T.
$$
We have
$$
|\Delta(t)| \le \frac{(t - T)^{\alpha +1}}{t^{\alpha +1}} \il{t_0}{T}|D(\tau)|d \tau + \frac{(\alpha +1)(t-t_0)^\alpha}{9 t^{\alpha +1}}\il{t_0}{T} p_{-}^2(\tau) d \tau + \phantom{aaaaaaaaaaaaaaaaaaaaa}
$$
$$
+\alpha(\alpha +1) c_T \frac{(t -T)^{\alpha -1}}{t^{\alpha +1}}\frac{(t - T)^2}{2} + \frac{\alpha(\alpha +1) y(T)(t-T)^{\alpha +1}}{t^{\alpha +1}}, \phh t \ge T. \eqno (3.6)
$$
Rewrite (3.5) in the form
$$
\frac{\alpha(\alpha +1)}{t^{\alpha +1}}\il{t_0}{t}(t - \tau)^{\alpha -1} y(\tau) d \tau \le - \frac{1}{t^{\alpha+1}}\il{t_0}{t}(t - \tau)^{\alpha}\biggl[(t -\tau) D(\tau) -  \frac{(\alpha +1)}{9} p_{-}^2(\tau)\biggr] d \tau + \phantom{aaaaaaaa}
$$
$$
\phantom{aaaaaaaaaaaaaaaaaaaaaaaaaaaaaaaaaaaaaaaaaaaaaaaaaaaa}+ \Delta(t), \phh t \ge T. \eqno (3.7)
$$
It follows from (3.5) that $\Delta(t)$ is a bounded function. Then from the condition $(b)$ of the theorem it follows that the right part of (3.7) takes negative values in some points $t \ge T$, whereas the hypothesis (3,2) implies that the left part of (3.2) is nonnegative for all $t \ge T.$ We have obtained a contradiction, completing the proof of the theorem.

Due to (3.1) from Theorem 3.1 we obtain immediately

{\bf Corollary 3.1.} {\it If $p(t) \equiv 0, \ph q(t) \le 0, \ph r(t) > 0, \ph t \ge t_0$ and for some $\alpha >1$
$$
\limsup\limits_{t \to +\infty}\frac{1}{t^\alpha}\il{t_0}{t}(t - \tau)^\alpha \Bigl[r(\tau) -\frac{2}{3\sqrt{3}} (-q(\tau))^{3/2}\Bigr]d \tau = +\infty,
$$
then Eq. (1.1) has oscillatory solutions.
}

{\bf Theorem 3.2.} {\it Let the condition $(a)$ of Theorem 3.1 and the following conditions be satisfied

\noindent
$ (c) \ph \ilp{t_0} D(\tau) d \tau = +\infty,$

\noindent
$(d)$ there exists $\alpha > 0$ such that
$$
\limsup\limits_{t \to+\infty}\frac{1}{t^{\alpha + 1}}\il{t_0}{t}(t - \tau)^{\alpha}\Bigl[(t - \tau) D(\tau) - (\alpha +1) p_{-}(\tau)\Bigr] d \tau = +\infty.
$$

\noindent
Then Eq. (1.1) has oscillatory solutions.
}

Proof. Let the condition $(a)$ of Theorem 3.1 holds. To prove the theorem it is enough, as in the proof of Theorem 3.1, to show that (3.2) cannot be satisfied.    Assume the contrary, that (3.2) is true. Then (3.3) holds. We set
$$
\lambda(t) \equiv - \il{T}{t}D(\tau) d \tau + c_T - y'(t) - \frac{3}{2}y^2(t) - p_{-}(t) y(t), \phh t \ge T.
$$
It follows from (3.3) that
$$
\lambda(t) \ge 0, \ph t \ge T. \eqno (3.8)
$$
Consider the Riccati equations
$$
y'+ \frac{3}{2}y^2 + p_{-}(t) y + \lambda(t) = - \il{T}{t} D(\tau) d \tau + c_T, \phh t \ge T, \eqno (3.9)
$$
$$
u' + \frac{3}{2} u^2 + p_{-}(t) u = 0, \phh t \ge T. \eqno (3.10)
$$
It follows from the condition $(c)$ that  there exists $T_1 \ge T$ such that
$$
\il{T}{t} D(\tau) d \tau \ge c_T, \phh t \ge T_1. \eqno (3.11)
$$
Obviously $y(t)$ is a solution of Eq. (3.9) on $[T_1,+\infty)$. Then applying Theorem 2.1 to the pair of equations (3.9) and (3.10) and taking into account (3.8) and (3.11) we conclude that the solution $u(t)$ of Eq. (3.10) with $u(T_1) > y(T_1)$ exists on  $[T_1,+\infty)$ and
$$
u(t) \ge y(t), \ph t \ge T_1. \eqno (3.12)
$$
It is well known that $u(t)$ can be represent in the following explicit form (as a solution of a Bernoully equation)
$$
u(t) = \frac{u(T_1)\exp\biggl\{- \il{T_1}{t} p_{-}(\tau) d \tau\biggr\}}{1 + \frac{3}{2}u(T_1)\il{T_1}{t}\exp\biggl\{-\il{T_1}{\tau} p_{-}(s) d s\biggr\} d \tau}, \phh t \ge T_1.
$$
Then from (3.12) it follows
$$
p_{-}(t) y(t) - \frac{p_{-}(t) u(T_1)\exp\biggl\{-\il{T_1}{t} p_{-}(\tau) d \tau\biggr\}}{1 + \frac{3}{2}u(T_1)\il{T_1}{t}\exp\biggl\{-\il{T_1}{\tau} p_{-}(s) d s\biggr\} d \tau} \ge 0, \phh t \ge T_1.
$$
This together with (3.3) implies
$$
y'(t) \le - \il{T}{t} D(\tau) d\tau + c_T -\frac{p_{-}(t) u(T_1)\exp\biggl\{-\il{T_1}{t} p_{-}(\tau) d \tau\biggr\}}{1 + \frac{3}{2}u(T_1)\il{T_1}{t}\exp\biggl\{-\il{T_1}{\tau} p_{-}(s) d s\biggr\} d \tau}.
$$
Integrating this inequality from $T_1$ to t we obtain
$$
y(t) \le - \il{T_1}{t} d\tau\il{T}{\tau} D(s) d s + \ln \biggl[1 + \frac{3}{2}u(T_1)\il{T_1}{t}\exp\biggl\{-\il{T_1}{\tau} p_{-}(s) d s\biggr\} d \tau\biggr] + y(T_1) + c_T(t - T),
$$
$t \ge T_1.$ From here it follows (since $p_{-}(t) \le 0, \ph t \ge t_0$)
$$
y(t) \le - \il{T_1}{t} d\tau\il{T}{\tau} D(s) d s + \ln \biggl[1 + \frac{3}{2}u(T_1)(t - T_1)\exp\biggl\{-\il{T_1}{t} p_{-}(s) d s\biggr\}\biggr] + y(T_1) + c_T(t - T),
$$
$t \ge T_1.$
For any $T_2 > T_1$ from here we obtain
$$
y(t) \le - \il{T_1}{t} d\tau\il{T}{\tau} D(s) d s  + \il{T_1}{t} p_{-}(s) d s   + L(t) +  \ln (t - T_1)  = c_T(t - T), \ph t \ge T_2, \eqno (3.13)
$$
where $L(t) \equiv \ln\biggl[1 + \frac{3}{2}u(T_1)(t - T_1)\exp\biggl\{-\il{T_1}{t} p_{-}(s) d s\biggr\}\biggr] - \il{T_1}{t} p_{-}(s) d s  - \ln (t - T_1) + y(T_1), \linebreak t \ge T_2.$ We claim that $L(t)$ is bounded on $[T_2,+\infty)$. Indeed, we have
$$
|L(t)| = \Biggl|\ln\biggl[1 + \frac{3}{2}u(T_1)(t - T_1)\exp\biggl\{-\il{T_1}{t} p_{-}(s) d s\biggr\}\biggr] -
$$
$$
-\ln\biggl[\frac{3}{2}u(T_1)(t - T_1)\exp\biggl\{-\il{T_1}{t} p_{-}(s) d s\biggr\}\biggr]  + \ln\Bigl[\frac{3}{2} u(T_1)\Bigr] + y(T_1)\Biggr| =
$$
$$
= \Biggl| \ln\Bigl[\frac{3}{2} u(T_1)\Bigr] + y(T_1)   + \ln \Biggl[1 + \frac{2}{3 u(T_1) \exp\biggl\{-\il{T_1}{t} p_{-}(s) d s\biggr\}(t- T_1)}\Biggr]\Biggr| \le
$$
$$
\le \biggl|\ln\biggl[\frac{3}{2}u(T_1)\biggr]\biggr| + y(T_1) + \ln \biggl[1  + \frac{2}{3 u(T_1) (T_2 - T_1)}\biggr] , \ph t \ge T_2.
$$
Therefore $L(t)$ is bounded on $[T_2,+\infty)$. From (3.13) we get
$$
y(t) \le - \il{T_2}{t} d\tau\il{T}{\tau} D(s) d s - \il{T_2}{t} p_{-}(s) d s + J(t), \phh t \ge T_2, \eqno (3.14)
$$
where
$$
J(t) \equiv L(t) -  \il{T_1}{T_2} d\tau\il{T}{\tau} D(s) d s - \il{T_1}{T_2} p_{-}(s) d s + \ln (t - T_1) - (t -T_2) \il{T}{T_2} D(s) d s + c_T (t - T), \ph t \ge T_2.
$$
Since $L(t)$ is bounded on $[T_2,+\infty)$ we have
$$
J(t) = O(t), \phh t \to +\infty. \eqno (3.15)
$$
Acting on both sides of (3.14) by the operator $K_{T_2.\alpha}$ we obtain
$$
\frac{\alpha(\alpha + 1)}{t^{\alpha + 1}}\il{T_2}{t}(t - \tau)^{\alpha -1} y(\tau) d \tau \le -\frac{1}{t^{\alpha + 1}}\il{T_2}{t}(t - \tau)^{\alpha + 1} D(\tau)d \tau +\phantom{aaaaaaaaaaaaaaaaaaaaaaaaaaa}
$$
$$
\phantom{aaaaaaaaaaaaaaa}+\frac{\alpha + 1}{t^{\alpha + 1}}\il{T_2}{t}(t - \tau)^{\alpha} p_{-}(\tau)d \tau + \frac{\alpha(\alpha + 1)}{t^{\alpha + 1}}\il{T_2}{t}(t - \tau)^{\alpha -1} J(\tau)d \tau, \phh t \ge T_2.
$$
Hence,
$$
\frac{\alpha(\alpha + 1)}{t^{\alpha + 1}}\il{T_2}{t}(t - \tau)^{\alpha -1} y(\tau) d \tau \le -\frac{1}{t^{\alpha + 1}}\il{t_0}{t}(t - \tau)^{\alpha + 1} D(\tau)d \tau +\phantom{aaaaaaaaaaaaaaaaaaaaaaaaaaa}
$$
$$
\phantom{aaaaaaaaaaaaaaaaaaa}+\frac{\alpha + 1}{t^{\alpha + 1}}\il{t_0}{t}(t - \tau)^{\alpha} p_{-}(\tau)d \tau + \Delta_1(t), \phh t \ge T_2, \eqno (3.16)
$$
where
$$
\Delta_1(t) \equiv \frac{1}{t^{\alpha + 1}}\il{t_0}{T_2}(t - \tau)^{\alpha + 1} D(\tau)d \tau -\phantom{aaaaaaaaaaaaaaaaaaaaaaaaaaaaaaaaaaaaaaaaaaaaa}
$$
$$
\phantom{aaaaaaaaaaaaaaa}-\frac{\alpha + 1}{t^{\alpha + 1}}\il{t_0}{T_2}(t - \tau)^{\alpha} p_{-}(\tau)d \tau - \frac{\alpha(\alpha + 1)}{t^{\alpha + 1}}\il{t_0}{T_2}(t - \tau)^{\alpha -1} J(\tau)d \tau, \phh t \ge T_2.
$$
It follows from (3.15) that $\Delta_1(t)$ is a bounded function on $[T_2,+\infty)$. Then it follows from the condition $(d)$ of the theorem that the right side of (3.16) takes negative values at some points $t$ of the interval $[T_2,+\infty)$, while the left side is always nonnegative. We have obtained a contradiction, which completes the proof of the theorem.

{\bf Example 3.1.} {\it Let $t_0 = 1, \ph p(t) = -M t^\gamma, \ph M > 0, \ph \gamma \in \mathbb{R}, \ph q(t) \equiv 0, \ph r(t) = N t^\beta - \delta(t), \ph  N > 0, \ph \beta > -1, \ph t \ge 1,$ where $\delta(t) = 0$, if $M t^{2\gamma} - 3\gamma t^{\gamma -1} < 0$ and $\delta(t) = \min\Bigl\{0, \Bigr[\frac{\sqrt{M^2 t^{2\gamma} - 3 M\gamma t^{\gamma - 1} + M t^\gamma}}{3}\Bigr]^3 - M t^\gamma\Bigr[\frac{\sqrt{M^2 t^{2\gamma} - 3 M\gamma t^{\gamma - 1} + M t^\gamma}}{3}\Bigr]^2 - M\gamma t^{\gamma -1} \frac{\sqrt{M^2 t^{2\gamma} - 3 M\gamma t^{\gamma - 1} + M t^\gamma}}{3}\Bigr\},$ if $M t^{2\gamma} - 3\gamma t^{\gamma -1} \ge 0, \ph t \ge 1.$ Then, obviously, $D(t) = N t^\beta, \ph p_{-}(t) = - M t^\gamma, \ph t \ge 1.$ Therefore, it is easy to check that, for $\gamma - 1 < \beta  < 2\gamma -1$ all the conditions of Theorem~ 3.2 are  satisfied, but the condition $(b)$ of Theorem 3.1 is not.
}

{\bf Theorem 3.3} {\it Let the condition $(a)$ of Theorem 3.1 and the following conditions be satisfied

\noindent
$(e) \ph \liminf\limits_{t \to +\infty} \il{t_0}{t} D(\tau) d\tau > -\infty,$

\noindent
$(f) \ph  p_{-}(t) = p_{-,1}(t) + p_{-,2}(t), \ph p_{-,j}(t)\le 0, \ph  j=1,2, \ph t \ge t_0, \ph \ilp{t_0} | p_{-,1}(t)| d t < +\infty$ and $p_{-,2}(t)$ is bounded,

\noindent
$(g)$ \ph for some $\alpha > 1$
$$
\limsup\limits_{t\to +\infty} \frac{1}{t^\alpha} \il{t_0}{t}(t -\tau)^\alpha D(\tau) d \tau = +\infty.
$$

\noindent
Then Eq. (1.1) has oscillatory  solutions.
}

Proof. By  the condition $(a)$ to prove the theorem it is enough to show that the relation (3.2) cannot be satisfied. Assume (3.2) is valid. Then (3.3) holds.
Rewrite (3.3) in the form
$$
y'(t) + \frac{3}{2}\Bigl(y(t) + \frac{p_{-,2}(t)}{3}\Bigr)^2 + p_{-,1}(t) y(t) \le -\il{T}{t}D(\tau) d \tau + c_T + \frac{p^2_{-,2}(t)}{9}, \ph t \ge T.
$$
It follows from here that
$$
y'(t) + p_{-,1}(t) y(t) \le -\il{T}{t}D(\tau) d \tau + c_T + \frac{p^2_{-,2}(t)}{9}, \ph t \ge T. \eqno (3.17)
$$
Since $p_{-,2}(t)$ is bounded it follows from the condition $(e)$ that
$$
\il{T}{t} D(\tau) d \tau - c_T - \frac{p_{-,2}(t)}{9} \ge - e_T \phh t \ge T, \eqno (3.18)
$$
for some $e_T > 0$. We set
$$
\lambda_1(t) \equiv -\il{T}{t} D(\tau) d \tau + c_T + \frac{p_{-,2}(t)}{9} -y'(t) - p_{-,1}(t) y(t), \ph t \ge T.
$$
It follows from (3.17) that
$$
\lambda_1(t) \ge 0, \phh t \ge T. \eqno(3.19)
$$
Consider the linear equations
$$
u'  + p_{-.1}(t) u   +  \lambda_1(t) +  \il{T}{t} D(\tau) d \tau - c_T - \frac{p_{-,2}(t)}{9}  = 0, \phh t \ge T. \eqno (3.20)
$$
$$
v'  + p_{-.1}(t) v - e_T = 0, \phh t \ge T. \eqno (3.21)
$$
Obviously $y(t)$ is a solution of Eq. (3.20) on $[T,+\infty)$ and
$$
v(t) \equiv \exp\biggl\{-\il{T}{t}p_{-,1}(\tau) d\tau\biggr\}\biggl[y(T) + e_T\il{T}{t}\exp\biggl\{\il{T}{\tau}p_{-,1}(s) ds\biggr\}d\tau\biggr], \phh t \ge T
$$
is a solution of Eq. (3.21) on $[T,+\infty)$ with $u(T) = y(T)$. Obviously we can interpret the equations (3.20) and (3.21) as Riccati equations with $\equiv 0$ coefficients of $u^2$ and $v^2$. Then applying Theorem 2.1 to these equations and taking into account (3.18) and (3.19) we conclude that $y(t) \le v(t), \ph t \ge T$. Hence,
$$
p_{-,1}(t) y(t) \ge p_{-,1}(t) v(t), \ph t \ge T  \ph (\mbox{since} \ph y(t) > 0, \ph p_{-,1}(t) \le 0 , \ph t \ge T).
$$
This together with (3.17) implies
$$
y'(t) \le - \il{T}{t}D(\tau) d \tau + f_T(t), \phh t \ge T,
$$
where
$$
f_T(t) \equiv c_T + \frac{p^2_{-,2}(t)}{9} - p_{-,1}(t) v(t), \phh t \ge T.
$$
Integrating the obtained inequality from $T$ to $t$ we obtain
$$
y(t) \le y(T) - \il{Y}{t}d\tau\il{T}{\tau}D(s) d s + \il{T}{t}f_T(\tau) d \tau, \phh t \ge T. \eqno (3.22)
$$
It is not difficult to verify  (by integrating by parts) that
$$
\il{T}{t}p_{-,1}(\tau) v(\tau) d \tau = \exp\biggl\{-\il{T}{t}p_{-,1}(\tau) d \tau\biggr\}\biggl[y(T) + e_T\il{T}{t}\exp\biggl\{\il{T}{\tau}p_{-,1}(s) d s\biggr\}\biggr]- \phantom{aaaaaaaaaaa}
$$
$$
 \phantom{aaaaaaaaaaaaaaaaaaaaaaaaaaaaaaaaaaaaaaaaaaaaaaaaaaaa}- y(T) - e_T(t-T), \ph t \ge T.
$$
Then, obviously,
$$
f_T(t) = O(t), \ph \mbox{as} \ph t \to +\infty. \eqno (3.23)
$$
Without loss of generality we can take that $T > 0$. Then acting on both sides of (3.22) by the operator $K_{T,\alpha-1}$ (here $\alpha >1$) we obtain
$$
\frac{(\alpha -1)\alpha}{t^\alpha}\il{T}{t}(t - \tau)^{\alpha -1} y(\tau) d \tau \le \frac{1}{t^\alpha}\il{T}{t}(t-\tau)^\alpha D(\tau) d\tau + \frac{(\alpha -1)\alpha}{t^\alpha}\il{T}{t}(t - \tau)^{\alpha -2}d\tau \il{T}{\tau} f_T(s) d s,
$$
$t \ge T$. Hence,
$$
\frac{(\alpha -1)\alpha}{t^\alpha}\il{T}{t}(t - \tau)^{\alpha -1} y(\tau) d \tau \le \frac{1}{t^\alpha}\il{t_0}{t}(t-\tau)^\alpha D(\tau) d\tau + F_T(t), \ph t \ge T, \eqno (3.24)
$$
where $F_T(t) \equiv -\frac{1}{t^\alpha}\il{t_0}{T}(t-\tau)^\alpha D(\tau) d\tau + \frac{(\alpha -1)\alpha}{t^\alpha}\il{T}{t}(t - \tau)^{\alpha -2} d\tau \il{T}{\tau}f_T(s) d s, \ph t \ge T.$ It follows from (3,23) that the function $F_T(t)$ is bounded on $[T,+\infty)$.  Obviously the left part of (3.24) is nonnegative for all $t \ge T$. However, it follows from the condition $(g)$ of the theorem that the right part of (3.24) takes negative values in some points of $[T,+\infty)$. We have obtained a contradiction, completing the proof of the theorem.

{\bf Remark 3.1.} {\it Theorems 3.1 - 3.3 are generalizations of Theorem 1.1.}

{\bf Example 3.2.} {\it Assume
$$
p(t) = M \sist{- n^3\sin^2[n^5\pi(t -n)], \ph t \in [n,n+\frac{1}{n^5}].}{0, \ph t \in (n+ \frac{1}{n^5}),}
$$
$M = const > 0, \ph  n= 1, 2, \dots $. Obviously $p(t)$ is continuously differentiable on $[1,+\infty)$ and $p_{-}(t) = p(t), \ph t\ge 1$. Let $q(t) \equiv 0, \ph r(t) = r_0 + \delta_1(t), \ph t \ge 1$, where $r_0 > 0, \ph \delta_1(t) = 0,$ if $p^2(t) + 3 p'(t) < 0$ and $\delta_1(t) = \min\Bigl\{0, \Bigl[\frac{\sqrt{p^2(t) + 3 p'(t)} - p(t)}{3}\Bigr]^3 + p(t)\Bigl[\frac{\sqrt{p^2(t) + 3 p'(t)} - p(t)}{3}\Bigr]^2 - p-(t) \frac{\sqrt{p^2(t) + 3 p'(t)} - p(t)}{3}\Bigr\},$ if $p^2(t) + 3 p'(t) \ge 0, \ph t \ge 1.$ Then obviously
$$
D(t) \equiv r_0 > 0, \phh t \ge1 \eqno (3.25)
$$
and $\ilp{t_0}|p_{-,1}(t)| d t = \ilp{t_0}|p(t)| d t  \le \sum\limits_{n=1}^{+\infty} \frac{1}{n^2} < +\infty$. Therefore all conditions of Theorem~ 3.3 for this case of coefficients  of Eq (1.1) are satisfied. Let us show that for this case the condition $(b)$ of Theorem 3.1 is not satisfied for all  $M \ge 12 r_0$. We have
$$
\il{n}{n+1}p_{-}^2(t) d t = n^6 \il{n}{n+\frac{1}{n^5}}\sin ^4 (n^5\pi(t - n))d t  = \frac{M n}{\pi} \il{0}{\pi}\sin^4 \tau d \tau = \frac{3 M n }{8}, \phh n =1,2,\dots .
$$
 It follows from here that
$$
\il{1}{t} p_{-}^2(\tau) d \tau \ge \frac{3 M(t -1) t}{16}, \phh t \ge 1.
$$
This together with (3.25) implies
$$
\il{1}{t}(t - \tau)^\alpha[(t -\tau) D(\tau) - \frac{\alpha +1}{9} p_{-}^2(\tau)]d \tau \le \il{1}{t}\frac{(t - \tau)^\alpha}{2}[r_0 \tau^2 - \frac{M(\alpha +1)}{12}  (\tau -1) \tau] d \tau, \ph t \ge 1.
$$
Therefore, if $M \ge 12 r_0$, then the condition $(b)$ of Theorem 3.1 is not satisfied,
}

\vskip 20pt

\centerline{\bf References}

\vskip 10pt

\noindent
1. L. Erbe, Existence of oscillatory solutions and asymptotic behavior for a class of third \linebreak \phantom{a} order linear differential equations. Pacific J. Math., vol. 64, No 2, 1976, pp. 369--385.

\noindent
2. G. A. Grigoriam, Some properties of the solutions of third order linear ordinary \linebreak \phantom{a} differential equations. Rocky J. Math. vol. 46, No 1, 2016, pp. 147--168.

\noindent
3. G. A. Grigorian,  On two comparison tests for second-order linear ordinary differential \linebreak \phantom{a}
equations, Diff. Urav., vol 47, (2011), 1225--1240 (in Russian), Diff. Eq., vol. 47, (2011),   \linebreak \phantom{a}
1237--1252  (in English).

\noindent
4. A. C. Lazer, The behavior of solutions of the differential equation $y''' + p(x)t' + q(x) y = 0$. \linebreak \phantom{a} Pacific J. Math., vol. 17, No 3, 1966, pp. 435--466.

\end{document}